\documentclass[12pt]{article}
\usepackage{amsfonts}
\usepackage{amssymb}
\usepackage{amsmath}
\usepackage{eufrak}
\usepackage[dvips]{graphicx}
\usepackage{latexsym}
\newcommand{\C}{{\EuFrak C}}

\newcommand{\bi}{\begin{itemize}}
\newcommand{\ei}{\end{itemize}}

\newtheorem{theorem}{Theorem}[section]

\newtheorem{fact}[theorem]{Fact}

\newtheorem{corollary}[theorem]{Corollary}
\newtheorem{proposition}[theorem]{Proposition}
\newtheorem{definition}[theorem]{Definition}
\newtheorem{remark}[theorem]{Remark}
\newtheorem{conjecture}[theorem]{Conjecture}
\newtheorem{question}[theorem]{Question}

\title{On $\omega$-categorical groups and rings with NIP}
\author{Krzysztof Krupi\'nski\footnote{Research supported by the Polish Government grant N N201 545938}}
\date{}

\begin{document}
\maketitle
\begin{abstract}
We show that $\omega$-categorical rings with NIP are nilpotent-by-finite.
We prove that an $\omega$-categorical group with NIP and fsg is nilpotent-by-finite.
We also notice that an $\omega$-categorical group with at least one strongly regular type is abelian.
Moreover, we get that each $\omega$-categorical, characteristically simple $p$-group with NIP has an infinite, definable abelian subgroup. Assuming additionally the existence of a non-algebraic, generically stable over $\emptyset$ type, such a group is abelian.
\end{abstract}

\section[\mbox{}]{Introduction}

Recall that a first order structure $M$ in a countable language is said to be $\omega$-categorical if, up to isomorphism, $Th(M)$ has at most one model of cardinality $\aleph_0$.
There is a long history of results describing the structure of $\omega$-categorical groups and rings. However, many questions in this area are still wide open. It follows easily that each countable, $\omega$-categorical group has a finite series of characteristic subgroups in which all successive quotients are characteristically simple groups. On the other hand, Wilson (see \cite{Wi,Ap}) proved 

\begin{fact}\label{Wilson}
For each countably infinite, $\omega$-categorical, characteristically simple group $H$, one of the following holds.
\begin{enumerate}
\item $H$ is an elementary abelian $p$-group for some prime $p$.
\item $H \cong B(F)$ or $H \cong B^{-1}(F)$ for some non-abelian, finite, simple group $F$, where $B(F)$ is the group of all continuous functions from the Cantor set $\cal{C}$ to $F$, and $B^-(F)$ is the subgroup of $B(F)$ consisting of the functions $f$ such that $f(x_0)=e$ for a fixed element $x_0 \in \cal{C}$.
\item $H$ is a perfect $p$-group.
\end{enumerate}
\end{fact}
Moreover, it is conjectured that (iii) is not realized.

As to $\omega$-categorical rings in general, we know that their Jacobson radical is nilpotent (see \cite[Lemma 1.3]{BR} and \cite{Ch}). However, there are examples of infinite, $\omega$-categorical rings which are semisimple (i.e. with trivial Jacobson radical) and so not nilpotent-by-finite \cite{BR}.

Interesting questions arise when one imposes additional model-theoretic restrictions (e.g. stability or simplicity) on our $\omega$-categorical group or ring. In the superstable [or, more generally, supersimple] $\omega$-categorical context, everything is clear: groups are [(finite central)-by-]abelian-by-finite \cite{Po,EW}; rigs are [(finite null)-by-]null-by-finite \cite{KW}. In the stable [or, more generally, NSOP] situation, we only know that $\omega$-categorical groups are nilpotent-by-finite \cite{Ma}, and $\omega$-categorical rings are nilpotent-by-finite \cite{BR, Kr1}, too. It is an open question whether $\omega$-categorical stable groups are abelian-by-finite and whether $\omega$-categorical stable rings are null-by-finite.

Our motivating problem is to describe the structure of $\omega$-categorical groups and rings satisfying NIP. Reasonable conjectures seems to be:

\begin{conjecture}\label{groups}
Each $\omega$-categorical group with NIP is nilpotent-by-finite.
\end{conjecture}

\begin{conjecture}\label{rings}
Each $\omega$-categorical ring with NIP is nilpotent-by-finite.
\end{conjecture} 

In this paper, we prove Conjecture \ref{rings}. As to Conjecture \ref{groups}, we prove it under the additional assumption that the group has fsg (finitely satisfiable generics). 

One of the ingredients of our proof of Conjecture \ref{groups} is the result saying that each $\omega$-categorical, characteristically simple $p$-group with NIP and having a non-algebraic, generically stable over $\emptyset$ type is abelian (so, under all these assumptions, (iii) of Fact \ref{Wilson} cannot happen).

We also show that each infinite, $\omega$-categorical, characteristically simple $p$-group with NIP has an infinite, definable abelian subgroup. Recall that the existence of an infinite, abelian subgroup is known for all infinite, $\omega$-categorical groups (use the fact that such groups are locally finite together with \cite[Corollary 2.5]{KW}). However, Plotkin \cite{Pl} found examples of  infinite, $\omega$-categorical $p$-groups with no infinite, definable, abelian subgroup.

At the end of the paper, we observe that an $\omega$-categorical group having at lest one strongly regular type is abelian.

I am grateful to Dugald Macpherson for interesting discussions and suggestions.

\section{Preliminaries}

Recall that we say that a group $G$ is solvable-by-finite [nilpotent-by-finite or abelian-by-finite] if it has a finite index (normal) subgroup which is solvable [nilpotent or abelian, respectively]. 
It is standard (see e.g. \cite[Remark 2.5]{Kr1})  that if $G$ is nilpotent-by-finite [abelian-by-finite], then it has a definable normal subgroup of finite index which is nilpotent [abelian, respectively]. If $G$ is solvable-by-finite, it is not clear whether it has a definable, solvable subgroup of finite index (it has such a subgroup if we additionally assume either icc on centralizers for all definable quotients of definable subgroups \cite[Remark 3.3]{Kr} or $\omega$-categoricity).

Recall some basic notions from ring theory. In this paper, rings are associative, but they are not assumed to contain 1  or to be commutative. An element $r$ of a ring $R$ is nilpotent of nilexponent $n$ if $r^n=0$, and $n$ is the smallest number with this
property. The ring is nil [of nilexponent $n$] if every element is nilpotent [of nilexponent $\leq n$, and there is an element of nilexponent $n$]. The  ring is nilpotent of class $n$ if $r_1\cdots r_n=0$ for all $r_1,\ldots,r_n\in R$, and $n$ is the smallest number with this property. An element $r$ is null if $rR=Rr=\{0\}$. The ring is null if all its elements are. 

We say that a ring $R$ is nilpotent-by-finite [null-by-finite] if it has a finite index ideal (equivalently subring by \cite{Le}) which is nilpotent [null, respectively]. 
By virtue of \cite[Remark 2.7]{Kr1}, this ideal can be chosen definable.

The Jacobson radical of a ring $R$ is the collection of all elements of $R$ satisfying the formula $\varphi(x)=\forall y \exists z (yx+z+zyx=0)$. The ring $R$ is semisimple if $J(R)=\{0\}$. It is always the case that $R/J(R)$ is semisimple.

Recall that a ring $R$ is a subdirect product of rings $R_i$, $i \in I$, if there is a monomorphism of $R$ into $\prod_{i \in I}R_i$ whose image projects onto each $R_i$. The following is \cite[Corollary 1]{BR}. 

\begin{fact}\label{fact1}
If $R$ is a semisimple, $\omega$-categorical ring, then $R$ is a subdirect product of complete matrix rings over finite fields. Moreover, only finitely many different matrix rings occur as subdirect factors.
\end{fact}

By \cite[Lemma 1.3]{BR} and \cite{Ch} we have
\begin{fact}\label{cherlin}
If $R$ is an $\omega$-categorical ring, then $J(R)$ is nilpotent.
\end{fact}

So, in order to prove that an $\omega$-categorical ring $R$ satisfying some extra assumptions is nilpotent-by-finite, it is enough to show that the semisimple ring $R/J(R)$ is finite (here Fact \ref{fact1} may be very helpful). We will use this approach in the proof of Theorem \ref{rings are nilpotent}.

We will also use \cite[Theorem 3.15]{Kr1}.

\begin{fact}\label{rings vs groups}
Suppose $G$ is a solvable, $\omega$-categorical group such that each ring interpretable in it is nilpotent-by-finite and each group $H$ definable in it has a definable connected component $H^0$ (i.e. the smallest definable subgroup of finite index). Then $G$ is nilpotent-by-finite. 
\end{fact}

Now, we recall the relevant notions from model theory. Let $T$ be a first order theory. We work in a monster model $\C$ of $T$.

\begin{definition}
We say that $T$ has the NIP if there is no formula $\varphi(x,y)$ and sequence $(a_i)_{i < \omega}$ such that for every $w \subseteq \omega$, there is $b_w$ such that $\models \varphi(a_i,b_w)$ iff $i \in w$.
\end{definition}

The next fact is Theorem 1.0.5 of \cite{Wa}.

\begin{fact}\label{wagner}
If $G$ is a group defined in a theory with NIP, then for each $\varphi$, there is some $N$ such that the intersection of any finite family of $\varphi$-definable subgroups of $G$ is an intersection of at most $N$ members of the family. 
\end{fact}

Let $p \in S(\C)$ be invariant over $A\subset \C$. We say that $(a_i)_{i \in \omega}$ is a Morley sequence in $p$ over $A$ if $a_n \models p |A a_{<n}$ for all $n$. It turns out that  Morley sequences  in $p$ over $A$ are indiscernible over $A$ and they have the same type over $A$.
If $\C' \succ \C$ is a bigger monster model, then the defining scheme of $p$ determines a unique $A$-invariant extension $\widetilde{p}\in S(\C')$ of $p$. By a Morley sequence in $p$ we mean a Morley sequence in $\widetilde{p}$ over $\C$. Finally, $p^{(k)}$ (where $k \in \omega \cup\{ \omega \}$) denotes the type over $\C$ of a Morley sequence in $p$ of length $k$

Recall from \cite{PT} that a global type $p \in S(\C)$ is said to be generically stable if, for some small $A$, it is $A$-invariant and for each Morley sequence $(a_i: i \in \omega)$ in $p$ over $A$ and each formula $\varphi(x)$ (with parameters from $\C$), the set $\{i : \models \varphi(a_i)\}$ is either finite or co-finite. This definition does not depend on the choice of $A$ over which $p$ is invariant. We will say that $p$ is generically stable over $A$ to express that $p$ is invariant over $A$ and generically stable. Assuming NIP, there are various equivalent definitions of generic-stability (see \cite[Proposition 3.2]{HP}). For us, one of them will be particularly important.

\begin{fact}\label{indiscernible set}
Assume $T$ has NIP, and $p \in S(\C)$ is $A$-invariant. Then, $p$ is generically stable iff every/some Morley sequence in $p$ over $A$ is an indiscernible set over $A$. 
\end{fact}

Now, we will briefly discuss fsg and generic stability. For more details on these and related notions see \cite{HPP,HP}. 

\begin{definition}
Let $G$ be a group definable in $\C$ by a formula $G(x)$.\\
(i) $G$ has fsg (finitely satisfiable generics) if there is a global type $p$ containing $G(x)$ and a model $M\prec \C$, of cardinality less than the degree of saturation of $\mathfrak C$, such that for all $g$, $gp$ is finitely satisfiable in $M$.\\
(ii) $G$ is generically stable if $G$ has fsg and some global generic type of $G$ is generically stable.
\end{definition}

We say that a group definable in a non-saturated model has one of the above properties if the group defined by the same formula in a monster model has it.

In general. generic stability of $G$ is a strictly stronger notion than fsg, but it is easy to check that these notions agree when $G^{00}$ (the smallest type-definable subgroup of bounded index) is definable and $T$ has NIP (see \cite[Section 6]{HP}). In the $\omega$-categorical, NIP context, $G^{00}$ is definable, and thus we get

\begin{remark}\label{fsg=gen stab}
Assume $G$ is a group definable in an $\omega$-categorical structure with NIP. Then, it has fsg iff it is generically stable.
\end{remark}

At the end, we recall the notion of a strongly regular type from \cite{PT}. 
The geometric meaning of this notion is explained in \cite{PT}.

\begin{definition}
Let $p(x) \in S(\C)$ be non-algebraic. We say that $p(x)$ is strongly regular if, for some small $A$, it is $A$-invariant and for all $B\supseteq A$ and $a$ from the sort of $x$, either $a \models p|B$ or $p|B \vdash p|Ba$.
\end{definition}

\section{$\boldsymbol{\omega}$-categorical rings with NIP}

In this section, we prove Conjecture \ref{rings}.

\begin{theorem}\label{rings are nilpotent}
Each $\omega$-categorical ring with NIP is nilpotent-by-finite.
\end{theorem}

\noindent
{\em Proof.}  By Fact \ref{cherlin}, everything boils down to showing that a semisimple, $\omega$-categorical ring $R$ with NIP is finite. Suppose for a contradiction that $R$ is infinite.

By $\omega$-categoricity, the two sided ideals $RrR$, $r \in R$, are uniformly definable (because $\omega$-categoricity implies that there is $K$ such that every element of any $RrR$ is the sum of at most $K$ elements of the form $r_1rr_2$ for $r_1,r_2 \in R$). Thus, by NIP and Fact \ref{wagner}, there is $N\geq 1$ such that for all $n \in \omega$ and  $r_0,\dots,r_n \in R$, there are $i_1,\dots, i_N \in \{0,\dots,n\}$ such that $Rr_0R\cap \dots \cap Rr_nR = Rr_{i_1}R \cap \dots \cap Rr_{i_N}R$.


Fact \ref{fact1} tells us that $R$ can be treated as a subring of the product $\prod_{i \in I}R_i$ of finite rings $R_i$ with identity, which projects onto each $R_i$, and where there are only finitely many pairwise distinct rings among the $R_i$'s, $i \in I$. Let $\pi_i$ be the projection onto the $i$th coordinate. 
For $i_0,\dots,i_n \in I$ and $r_j \in R_{i_j}$, we introduce  the set 
$$R_{i_0,\dots,i_n}^{r_0,\dots,r_n}=\left\{ r \in R: \bigwedge_{j=0}^n \pi_{i_j}(r)=r_{j}\right\}.$$

Using the assumption that $R$ is infinite and $R_i$'s are finite, we see that for any $i_0,\dots, i_n \in I$,
\begin{equation*}\tag{$*$}
R_{i_0,\dots,i_n}^{0,\dots,0}\;\, \mbox{is an infinite ideal of R.}
\end{equation*}
{\bf Claim 1}
There are  pairwise distinct $i_0, i_1,\dots \in I$, non-nilpotent elements $r_j\in R_{i_j}$, and elements $\eta_j \in R$ such that for all $n \in \omega$,
$$\eta_n \in R_{i_0,\dots,i_n}^{0,\dots,0,r_n}.$$
{\em Proof of Claim 1.} Suppose $i_0,\dots, i_{n}$, $r_0,\dots,r_n$ and $\eta_0,\dots,\eta_n$ have been chosen. Since $R$ is semisimple, it has no non-trivial nil ideals. Thus, by (*), $ R_{i_0,\dots,i_n}^{0,\dots,0}$ contains a non-nilpotent element $\eta_{n+1}$. As there are only finitely many different $R_i$'s and they are all finite, we can find $i_{n+1}$ different from $i_0,\dots,i_n$ such that $r_{n+1}:=\pi_{i_{n+1}}(\eta_{n+1}) \in R_{i_{n+1}}$ is non-nilpotent. \hfill $\blacksquare$\\  

%
\noindent
{\bf Claim 2} There are natural numbers $n(0)<\dots<n(N)$ such that the sets $$R_{i_{n(0)},\dots,i_{n(N)}}^{r_{n(0)},0,\dots,0}, R_{i_{n(0)},\dots,i_{n(N)}}^{0,r_{n(1)},0,\dots,0}, \dots, R_{i_{n(0)},\dots,i_{n(N)}}^{0,\dots,0, r_{n(N)}}$$ are non-empty.\\

Before proving Claim 2, let us notice that it leads to a contradiction. Choose $a_0,\dots,a_N$ from $R_{i_{n(0)},\dots,i_{n(N)}}^{r_{n(0)},0,\dots,0},\dots, R_{i_{n(0)},\dots,i_{n(N)}}^{0,\dots,0, r_{n(N)}}$, respectively. Put $b_k = \sum_{l \ne k}a_l$ for $k=0,\dots,N$. Then, 
\begin{equation*}\tag{$**$}
\pi_{i_{n(j)}}[Rb_0R\cap\dots \cap Rb_NR]= \{0\}\;\,  \mbox{for}\;\, j=0,\dots, N. 
\end{equation*}
On the other hand, $\prod_{k \ne j}b_k \in \bigcap_{k \ne j}Rb_kR$ for $j=0,\dots,N$. We also have that $\pi_{i_{n(j)}}[\prod_{k \ne j}b_k]=r_{n(j)}^N \ne 0$ as $r_{n(j)}$ is non-nilpotent. So, 
\begin{equation*}\tag{$***$}
 \pi_{i_{n(j)}}[\bigcap_{k \ne j}Rb_kR] \ne \{ 0 \} \;\,  \mbox{for}\;\, j=0,\dots, N.
\end{equation*}
 By $(**)$ and $(***)$,  $Rb_0R\cap \dots\cap Rb_N R \ne \bigcap_{k\ne j}Rb_kR$  for all $j=0,\dots,N$. This is a contradiction with the choice of $N$.\\


\noindent
{\em Proof of Claim 2.} 
%
%
%
%
%
Let $c=\max_{i \in I} |R_i|$. Define recursively:
$$
\begin{array}{lll}
L_N & = & c+1,\\
L_{N-k} & = & c^{L_N+\dots+L_{N-k+1}+1}+1 \;\mbox{for}\; k=1,\dots,N-1.
\end{array}
$$
Put $L_0=0$.

Consider any $k\in \{0,\dots,N-1\}$. Suppose each natural number $\alpha$ from the closed interval $[L_{N-k-1}+\dots + L_0, L_{N-k}+\dots +L_0 -1]$ has color
$$(\pi_{i_{L_{N-k}+\dots+L_1}}(\eta_\alpha),\dots,\pi_{i_{L_{N}+\dots+L_1}}(\eta_\alpha))\in \prod_{j=L_{N-k}+\dots+L_1}^{L_N+\dots+L_1} R_{i_j}.$$

In this way, the natural numbers from the interval $[L_{N-k-1}+\dots + L_0, L_{N-k}+\dots +L_0 -1]$ have been colored with at most $c^{L_N+\dots+ L_{N-k+1}+1}=L_{N-k}-1$ colors. Since there are $L_{N-k}$ such numbers, we can find natural numbers $n(N-k-1)<n'(N-k-1)$ from $[L_{N-k-1}+\dots + L_0, L_{N-k}+\dots +L_0 -1]$ with the same color. Put $a_{N-k-1}= \eta_{n(N-k-1)} -\eta_{n'(N-k-1)}$. Then, $\pi_{i_{j}}(a_{N-k-1})=0$ for all $j \in [L_N+\dots+L_{N-k}, L_N+\dots +L_1]$. Moreover, from the choice of $\eta_n$'s, we see that  $\pi_{i_{j}}(a_{N-k-1})=0$ for all $j <L_{N-k-1}+\dots +L_0$.
Putting additionally $n(N)=L_N+\dots+L_1$, we get a sequence $n(0)<\dots<n(N)$ which satisfies the conclusion of the claim. \hfill $\blacksquare$

\section{$\boldsymbol{\omega}$-categorical groups with NIP}

In this section, we investigate the structure of $\omega$-categorical groups with NIP. First, we make some observations on characteristically simple groups in this context. Then, we prove Conjecture \ref{groups} under the additional assumptions of fsg or of the existence of a strongly regular type.

As was mentioned in the introduction, each $\omega$-categorical group is locally finite, and so, if it is infinite, it has an infinite, abelian subgroup \cite[Corollary 2.5]{KW}. 
However, \cite{Pl} yields examples of infinite, $\omega$-categorical $p$-groups with no infinite, definable, abelian subgroup. 

\begin{proposition}
Let $p$ be a prime number. Then every infinite, characteristically simple, $\omega$-categorical $p$-group $G$ with NIP has an infinite, definable, abelian subgroup.
\end{proposition}
{\em Proof.} For any $x_0,\dots,x_n \in G$, $\langle x_0,\dots,x_n \rangle$ is a finite $p$-group, and so $C(x_0)\cap \dots \cap C(x_n) \ne \{ e \}$. Hence, if there were $x_0,\dots,x_n \in G$ with $C(x_0)\cap \dots \cap C(x_n)$ finite, we would find $x_{n+1},\dots,x_m \in G$ such that $Z(G)=C(x_0)\cap \dots \cap C(x_m) \ne \{ e \}$, which would  contradict the characteristic simplicity of $G$. So, we have proved the following\\[3mm]
{\bf Claim} For any $n \in \omega$ and $x_0,\dots,x_n \in G$, the intersection $C(x_0)\cap \dots \cap C(x_n)$ is infinite.\\

By the Claim, we can choose a sequence $(x_n)_{n\in \omega}$ of pairwise distinct elements of $G$ such that $x_{n+1} \in C(x_0)\cap \dots \cap C(x_n)$ for all $n \in \omega$. Put $G_n = C(x_0) \cap \dots \cap C(x_n)$. Then, $x_0, \dots, x_n \in Z(G_n)$. Hence, $|Z(G_n)|>n$.
On the other hand, by NIP and Fact \ref{wagner}, there is $N$ such that for any $y_0,\dots,y_n \in G$, there are $i_1,\dots, i_N \in \{0,\dots,n\}$ with $C(y_0)\cap \dots \cap C(y_n)=C(y_{i_1}) \cap \dots \cap C(y_{i_N})$. Thus, $G_n = C(x_{i_1^n}) \cap \dots \cap C(x_{i_N^n})$ for some $i_1^n, \dots, i_N^n \in \{ 0,\dots, n\}$. Since by $\omega$-categoricity the set $\{ tp(x_{i_1^n},\dots,x_{i_N^n}): n \in \omega\}$ is finite, there is $n \in \omega$ such that $Z(G_n)$ is infinite. \hfill $\blacksquare$\\

The next proposition uses notation from Fact \ref{Wilson}.

\begin{proposition}\label{B(F)}
For any non-abelian, finite, simple group $F$, neither $B(F)$ nor $B^-(F)$ have NIP.
\end{proposition}
{\em Proof.} 
Let $C_0,C_1,\dots$ be disjoint clopen subsets of the Cantor set $\cal{C}$ not containing $x_0$. Choose $g \in F \setminus Z(F)$. Define a sequence $(f_i)_{i \in \omega}$ of elements of $B^-(F)$ by:
$$
f_i(\eta)=\left\{ \begin{array}{ll}
g & \mbox{if}\; \eta \in C_i,\\
e & \mbox{if}\;  \eta \notin C_i.
\end{array} \right.
$$

Now, suppose for a contradiction that $B(F)$ has NIP (the case when $B^-(F)$ has NIP is almost the same). Using NIP and Fact \ref{wagner}, and reordering $C_i$'s if necessary, we can find $N$ such that $C_{B(F)}(f_0) \cap \dots \cap C_{B(F)}(f_N)= C_{B(F)}(f_0) \cap \dots \cap C_{B(F)}(f_{N-1})$. 

Take $h \in F \setminus C(g)$ and define $f \in B^-(F)$ by:
$$
f(\eta)=\left\{ \begin{array}{ll}
h & \mbox{if}\; \eta \in C_N,\\
e & \mbox{if}\;  \eta \notin C_N.
\end{array} \right.
$$
Then, we see that $f \in C_{B(F)}(f_0) \cap \dots \cap C_{B(F)}(f_{N-1}) \setminus C_{B(F)}(f_0) \cap \dots \cap C_{B(F)}(f_{N})$, a contradiction. \hfill $\blacksquare$

\begin{proposition}\label{p-group}
Let $p$ be a prime number. Let $G$ be a characteristically simple $p$-group interpretable over $\emptyset$ in an $\omega$-categorical structure with NIP. Assume that $G$ has a global, non-algebraic type $q$ which is generically stable over $\emptyset$. Then $G$ is abelian. 
\end{proposition}
{\em Proof.}
Assume for simplicity that $G$ is the ambient structure.
Wlog $G$ is a monster model.
Let $(a_i)_{i \in \omega}$ be a Morley sequence in $q$ over $\emptyset$. By NIP and Fact \ref{wagner}, there is $N$ such that for any $m$, $C(a_0)\cap \dots \cap C(a_m)=C(a_{i_1})\cap \dots \cap C(a_{i_N})$ for some $i_1,\dots,i_N \in \{ 0, \dots,m\}$. But, using Fact \ref{indiscernible set}, $(a_i)_{i \in \omega}$ is an indiscernible set over $\emptyset$. This implies that for any $m \geq N-1$ and $0\leq i_1<\dots <i_N \leq m$, one has $C(a_0)\cap \dots \cap C(a_m)=C(a_{i_1})\cap \dots \cap C(a_{i_N})$.

Consider any $(b_0,b_1,\dots) \models q^{(\omega)} |\emptyset$. Then, there is a sequence $(c_i)_{i \in \omega}$ such that the sequences $(a_i : i \in \omega)^\frown (c_i : i \in \omega)$ and $(b_i : i \in \omega)^\frown (c_i : i \in \omega)$ are indiscernible over $\emptyset$. Thus, by the last paragraph, $\bigcap_{i \in \omega} C(a_i)=C(c_0)\cap \dots \cap C(c_{N-1})=\bigcap_{i \in \omega}C(b_i)$. But, $\langle c_0,\dots,c_{N-1} \rangle$ is a finite $p$-group, which implies that $C(c_0)\cap \dots \cap C(c_{N-1}) \ne \{ e \}$. We conclude that $\bigcap_{i \in \omega} C(a_i)$ is a non-trivial, $\emptyset$-invariant (so $\emptyset$-definable) subgroup. Since $G$ is characteristically simple, we get $\bigcap_{i \in \omega} C(a_i)=G$, which implies that $Z(G) \ne \{ e \}$, and so $G=Z(G)$ once again by the characteristic simplicity of $G$. \hfill $\blacksquare$

\begin{theorem}
Each $\omega$-categorical group $G$ with NIP and fsg is nilpotent-by-finite.
\end{theorem}
{\em Proof.}
We can assume that $G$ is infinite and it is a monster model. In fact, we will work in a slightly more general context, where $G$ is a group with fsg which is interpretable (over $\emptyset$)  in an $\omega$-categorical monster model $\C=\C^{eq}$ satisfying NIP.
Since $G^{00}$ is $\emptyset$-invariant, it is $\emptyset$-definable by $\omega$-categoricity. Thus, we can assume that $G=G^{00}$. Then, $G$ has a unique global generic type \cite[Proposition 0.26]{EKP}, which must be $\emptyset$-invariant.

By $\omega$-categoricty, there is a normal series $\{ e \}=G_0\le G_1\le \dots \le G_n=G$ of $\emptyset$-definable (in $\C$) subgroups of $G$ such that each $G_{i+1}/G_i$ is characteristically simple in the sense of $(G,\C)$, i.e. it does not have non-trivial, proper, subgroups normal in $G/G_i$ and invariant under $Aut(\C)$. In particular,  each $G_{i+1}/G_i$ is a characteristically simple group. Indeed, it is clear that any characteristic subgroup of $G_{i+1}/G_i$ is invariant under $Aut(\C)$; the fact that such a subgroup is normal in $G/G_i$ follows from the fact that inner automorphisms of $G$ induce automorphisms of $G_{i+1}/G_i$.   

Our proof will be by induction on $n$ (where $n$ is a number for which there exists a sequence as in the last paragraph).
Consider the case $n=1$. Then, $G=G_1$ is an infinite, characteristically simple group.
Proposition \ref{B(F)} eliminates the possibility that a countable elementary substructure of $G$ is as in point (ii) of Fact \ref{Wilson}. Proposition \ref{p-group} together with Remark \ref{fsg=gen stab} and the $\emptyset$-invariance of the unique generic type of $G$ eliminate the possibility from point (iii) of Fact \ref{Wilson}.  Thus, $G$ must be abelian.

We turn to the induction step, where we assume that $n\geq 2$.
First, notice that $\omega$-categoricity implies that the upper central series $Z_0(G)\le Z_1(G)\le Z_2(G)\le \dots$ stabilizes after finitely many, say $m$,  steps.
So, replacing $G$ by $G/Z_m(G)$, we can assume that $G$ is centerless. Such a replacement is possible, because the normal series $G_0Z_m(G)/Z_m(G)\le G_1Z_m(G)/Z_m(G)\le \dots\le G_nZ_m(G)/Z_m(G)=G/Z_m(G)$ consists of groups which are $\emptyset$-definable in $\C$, and it is easy to check that each quotient $(G_{i+1}Z_m(G)/Z_m(G))/(G_iZ_m(G)/Z_m(G))$ is characteristically simple in the sense of $(G/Z_m(G),\C)$.

Since $n \geq 2$, by the induction hypothesis, we can assume that $G_1$ is a non-trivial, proper subgroup of $G$. Since $G/G_1$ also has fsg, using the induction hypothesis, we get that $G/G_1$ is nilpotent-by-finite (notice that since we consider $\emptyset$-definability in $\C$, the normal series $G_1/G_1\le G_2/G_1\le \dots\le G_n/G_1=G/G_1$ allows us to use the induction hypothesis). So, it is nilpotent, because the Fitting subgroup of $G/G_1$ is a nilpotent subgroup of finite index which is $\emptyset$-definable and so equal to $G/G_1$ by the connectedness of $G/G_1$.

Let $q$ be the unique global generic type of $G$, and $(g_i)_{i \in \omega}$ be a Morley sequence in $q$ over $\emptyset$. Since $G/G_1$ is nilpotent, there is a minimal $k$ such that the iterated commutator $[g_{k-1},[g_{k-2},\dots,[g_1,g_0]\dots]] \in G_1$. Since $g_0$ is generic over $\emptyset$ and $G/G_1$ is infinite, we see that $k\geq 2$. Define 
$$h_i=[g_{ik+k-1},[g_{ik+k-2},\dots,[g_{ik+1},g_{ik}]\dots]]$$ 
for $i \in \omega$. 
Let $(g_i')_{ i \in \omega}$ be a Morley sequence in $q$ over $\C$. Put 
$$h_i'=[g_{ik+k-1}',[g_{ik+k-2}',\dots,[g_{ik+1}',g_{ik}']\dots]]$$ for $i \in \omega$. Since $tp(g_{k-1}',\dots,g_0'/\C)=q^{(k)}$ is invariant over $\emptyset$, the type $r:=tp(h_0'/\C)$ is also invariant over $\emptyset$. Moreover, $(h_i)_{i \in \omega}$ is a Morley sequence in $r$ over $\emptyset$. By the generic stability of $q$ and Fact \ref{indiscernible set}, the sequence $(g_i)_{i \in \omega}$ is an indiscernible set, and so $(h_i)_{i \in \omega}$ is an indiscernible set as well. We conclude that $r$ is generically stable over $\emptyset$. 

We claim that $r$ is non-algebraic. To see this, let us define 
$$a_j:=[g_{k-1+j}',[g_{k-2}',\dots,[g_1',g_0']\dots]]$$
for $j \in \omega$. We see that $a_j \models r$. Moreover, $a_{j_1} \ne a_{j_2}$ whenever $j_1<j_2$. Indeed, if this is not the case, then for some $j_1<j_2$, we have that $g_{j_2}'g_{j_1}'^{-1} \in C(h)$, where $h=[g_{k-2}',\dots,[g_1',g_0']\dots]$. By the minimality of $k$, $h \ne e$. But, $g_{j_2}'$ is generic over $g_{j_1}',h$, and so $g_{j_2}'g_{j_1}'^{-1}$ is generic over $h$. Thus, $[G(\C'):C(h)]<\omega$, where $\C'$ is a monster model containing $\C$ and $(g_i')_{i \in \omega}$. This implies $h \in Z(G(\C'))=\{ e \}$, a contradiction.

We have proved that $r$ is a global, non-algebraic type of $G_1$ which is generically stable over $\emptyset$ (in particular, $G_1$ is infinite). So, by Proposition \ref{p-group} together with Proposition \ref{B(F)} and Fact \ref{Wilson}, $G_1$ is abelian. Hence, $G$ is solvable. 

By Theorem \ref{rings are nilpotent} and Fact \ref{rings vs groups}, we conclude that $G$ is nilpotent-by-finite. \hfill $\blacksquare$\\

Dugald Macpherson told me an alternative ending of the above proof, i.e. an alternative proof of the fact that a solvable, $\omega$-categorical group with NIP is nilpotent-by-finite. Namely, by \cite{AM}, we know that each countable, solvable, $\omega$-categorical group which is not nilpotent-by-finite interprets the countable, atomless Boolean algebra. So, it remains to show that this algebra does not have NIP, which is an easy exercise.  

Now, we will drop the NIP and fsg assumption, and instead we will assume the existence of a strongly regular type. Recall the following question from \cite{PT}.

\begin{question}
Suppose $G$ is a group with at least one strongly regular type. Does it imply that $G$ is abelian?
\end{question}

\begin{proposition}
If $G$ is any group with at least one strongly regular type, then all non-central elements of $G$ are conjugated.
\end{proposition}
{\em Proof.} Taking an elementary extension of $G$, we can assume that there is a global type $p$ whose strong regularity is witnessed over $G$. Consider any non-central element $a \in G$. Take $b \models p|G$. 

Notice that if $a^b \models p|G$, then the formula defining the conjugacy class of $a$ belongs to $p|G$. Thus, all elements $a \in G$ for which $a^b \models p|G$ are in one conjugacy class.
So, it remains to show that the assumption $a^b \not\models p|G$ leads to contradiction.   

This assumption and the strog regularity of $p$ over $G$ imply that $tp(b/G) \vdash tp(b/a^b,G)$. Thus, there is a formula $\varphi(x,y)$ (without parameters) and $g \in G^n$ such that $b \models \varphi(x,g)$ and $\models (\varphi(x,g) \rightarrow a^x=a^b)$. So, there is $c \in G$ such that $a^c=a^b$, and hence $b \in C(a)c$.  This means that $p|G \vdash \,\mbox{'}x \in C(a)c\mbox{'}$.

Consider two distinct realizations $g_1$ and $g_2$ of $p|G$ (they exist because $p$ is non-algebraic).\\
{\bf Case 1} $c \in C(a)$.\\
Then, $p|G \vdash \, \mbox{'} x \in C(a) \mbox{'}$, so $g_1 \in C(a)$. Take $h \notin C(a)$. Then, $hg_1 \notin C(a)$. Thus,
by the strong regularity of $p$ over $G$, we conclude that $tp(g_1/G,h,hg_1) =p|G,h,hg_1=tp(g_2/G,h,hg_1)$. But, the formula $x=h^{-1}hg_1$ belongs to  $tp(g_1/G,h,hg_1)$ and does not belong to $tp(g_2/G,h,hg_1)$, a contradiction.\\
{\bf Case 2} $c \notin C(a)$.\\
Since $g_1 \in C(a)c$, we have $g_1c^{-1} \in C(a)$, and so $g_1c^{-1} \notin C(a)c$. Hence, by the strong regularity of $p$, we get $tp(g_1/G,g_1c^{-1})=p|G,g_1c^{-1}=tp(g_2/G,g_1c^{-1})$. But, the formula $x=g_1c^{-1}c$ belongs to $tp(g_1/G,g_1c^{-1})$ and does not belong to $tp(g_2/G,g_1c^{-1})$, a contradiction. \hfill $\blacksquare$

\begin{corollary}
If $G$ is a group with at least one strongly regular type, then all non-central elements of $G$ have infinite order. 
In particular, an $\omega$-categorical group with at least one strongly regular type is abelian.
\end{corollary}     
{\em Proof.} This is a standard argument. We can assume that $G\ne Z(G)$. Suppose for a contradiction that there is a non-central element of finite order. By the last proposition, $G/Z(G)$ has one non-trivial conjugacy class. So, all non-trivial elements of $G/Z(G)$ have the same order, which must be a prime number $p$. If $p=2$, then $G/Z(G)$ is abelian, so $[G:Z(G)]\leq 2$, which implies $G=Z(G)$, a contradiction. Now, we assume that $p$ is odd. Take  a nontrivial $g \in G/Z(G)$. Then, there is $h \in G/Z(G)$ such that $h^{-1}gh=g^{-1}$. So, $g \in C(h^2)\setminus C(h)$. Finally we get $$C(h)=C(h^{2^{p-1}}) \supsetneq C(h^{2^{p-2}}) \supsetneq  \dots \supsetneq C(h^2) \supsetneq C(h),$$ which is impossible.\hfill $\blacksquare$

\noindent
Krzysztof Krupi\'nski\\
Instytut Matematyczny Uniwersytetu Wroc\l awskiego\\
pl. Grunwaldzki 2/4, 50-384 Wroc\l aw, Poland.\\
e-mail: kkrup@math.uni.wroc.pl\\[5mm]

\end{document}